\documentclass[12pt]{amsart}
\usepackage{geometry} 
\usepackage{amsmath,amsthm,amsfonts,amssymb}
\usepackage{graphics,xy}
\usepackage{xcolor}
\usepackage{graphicx}
\usepackage[normalem]{ulem}
\usepackage{url}
\usepackage{tikz}
\usetikzlibrary{arrows.meta}

\tikzset{
	diagram/.style={line width=.9pt,line cap=round,line join=round},
	coupon/.style={draw,fill=white,rectangle,
		minimum width=6mm,minimum height=2.4mm,
		inner sep=0pt,outer sep=0pt}
}

\tikzset{
	diagram/.style={line width=.9pt,line cap=round,line join=round},
	spine/.style={draw=red},
	meridian/.style={draw=blue},
	hidden/.style={dash pattern=on 2.2pt off 2.2pt,line width=.65pt}
}

\usepackage{hyperref}
\usepackage[backend=biber, backref=true, maxbibnames=99, style=alphabetic]{biblatex}
\usepackage[utf8]{inputenc}

\newtheorem{lemma}{Lemma}[section]

\newtheorem{theorem}[lemma]{Theorem}

\theoremstyle{definition}
\newtheorem{definition}[lemma]{Definition}
\newtheorem{remark}[lemma]{Remark}

\newcommand{\Q}{\mathbb{Q}}
\newcommand{\C}{\mathbb{C}}
\newcommand{\N}{\mathbb{N}}

\newcommand{\ra}{\longrightarrow}

\title[Dimension bounds on skein modules]{Dimension bounds on skein modules}

\author{Giulio Belletti}
\address{
	IRMP, UC Louvain, Chemin du Cyclotron 2, bte L7.01.02, 1348 Louvain-la-Neuve, Belgium}
\email{gbelletti451@gmail.com}
\date{} 

\begin{document}
\maketitle
	\begin{abstract} 
		We show that for any closed manifold $M$ obtained via a genus $g\geq 2$ Heegaard splitting, the dimension of $S(M,\Q(A))$ is at most $B^{3g-3}$, where $B$ is a concrete count of intersection points coming from the Heegaard splitting. The proof is independent of the recent proofs of the Witten finiteness conjecture of \cite{BD25} and \cite{JR25}, and thus can be seen as a more streamlined verification of the conjecture. Furthermore, the proof also works when choosing as base ring $\C$ and as $A$ any non-zero complex number that is not a root of unity.
\end{abstract}
\tableofcontents

\section{Introduction}

The problem of whether the Kauffman bracket skein module $S(M,\Q(A))$ of a closed manifold $M$ is finite dimensional first appeared in print in \cite{Car}; it had originated from a question that Witten had asked some years prior and thus became known as Witten's finiteness conjecture. 

A related result in this direction first appeared in \cite{GJS19}, showing that a certain $h$-adic completion of the Kauffman bracket skein module is finite dimensional. Later, alternative proofs for the Kauffman bracket skein module were provided by \cite{BD25} and \cite{JR25}; the proofs were achieved independently and are very different from each other. Neither of these proofs provide a concrete estimate of the actual dimension of $S(M,\Q(A))$, although \cite{BD25} does provide an explicit algorithm to obtain a finite generating set.

The main result of this paper is the following dimension bound on the Kauffman bracket skein module of a closed oriented $3$-manifold:

\begin{theorem}\label{thm:main}
	Let $M$ be a closed oriented $3$-manifold, and $R$ either $\Q(A)$ or $\C$ with $A\in \C^*$ not a root of unity. Then for any Heegaard splitting $M=H\cup H'$ of genus $g\geq 2$, there exists a number $B\in \N$, explicitly computable from the Heegaard splitting, such that 
	$$\dim S(M,R)\leq B^{3g-3}.$$
\end{theorem}

\begin{remark}
	The proof of Theorem \ref{thm:main} does not rely on the previous finiteness results of \cite{BD25} or \cite{JR25}; it is therefore an alternative proof of the Witten finiteness conjecture, albeit following a similar proof strategy as \cite{BD25}. Furthermore, the statement works when $A$ is not a root of unity, extending the Witten finiteness conjecture to this case.
\end{remark}

\begin{remark}
	Remark 1.2 of \cite{Jor24} says it is expected that $S(M,R)$ is finite dimensional whenever $A$ is not a root of unity, and that the dimension should not depend on $A$. Theorem \ref{thm:main} confirms the first part of this expectation; the upper bound we establish does not depend on $A$ but this does not directly confirm that the dimension itself does not depend on $A$.
\end{remark}
Using the same techniques used to prove Theorem \ref{thm:main}, we can also obtain an upper bound on the dimension of the localized skein module for manifolds with boundary (\cite{BD25}); for simplicity, we focus on the closed case, and we will explain in Section \ref{sec:boundary} the minor modifications needed to deal with the boundary case.

\begin{remark}
	The proof of Theorem \ref{thm:main} follows a similar strategy to the proof in \cite{BD25}, although the execution of the strategy is much simpler. The proof of \cite{BD25}, however, is constructive: following the steps provides a finite generating set for $S(M,\Q(A))$, whereas the proof in this paper does not. 
\end{remark}

A forthcoming paper will use the techniques of this paper to prove finiteness of skein modules for higher rank Lie algebras.

\textbf{AI disclaimer:} The author used Astra 6 to assist with the proofs of this paper and with the generation of the TikZ images. Specifically, Astra provided the key idea to use the boundary operators to break the growth of the filtration of a handlebody, inverting the core strategy of \cite{BD25}, which was to use the growth of the filtration to produce boundary operators that act trivially. Apart from the figures, everything in the present manuscript was written by the author, and he takes full responsibility for the mathematical content of the paper. After writing, ChatGPT was used to check for mistakes, grammar, typos and so on.

\textbf{Acknowledgments}
 I wish to thank Renaud Detcherry, David Jordan, Thang Le, Iordanis Romaidis and Adam Sikora for comments and conversations about this work.
 
I am supported by FNRS (grant 1.B.044.25F). Part of this research was performed while visiting the Institute for Pure and Applied Mathematics (IPAM), which is supported by the National Science Foundation (Grant Nos. DMS-1925919 and DMS-2422832).
\section{Definitions and notation}
\begin{figure}
\begin{tikzpicture}[diagram,rotate=90,transform shape]
	\draw (-1.2,1.1) -- (0,0) -- (-1.2,-1.1);
	\draw (0,0) -- (1.3,0);
	
	\draw[-{Stealth[length=2.4mm]}] (0,-1.6) -- (0,-2.6);
	
	\begin{scope}[yshift=-4.2cm]
		\draw (-.83,.71)
		.. controls (-.36,.24) and (-.36,-.24) .. (-.83,-.71);
		\draw (-.61,.93)
		.. controls (-.25,.48) and (.24,.18) .. (.93,.18);
		\draw (-.61,-.93)
		.. controls (-.25,-.48) and (.24,-.18) .. (.93,-.18);
		
		\draw (-1.27,1.37) -- (-.72,.82);
		\draw (-1.27,-1.37) -- (-.72,-.82);
		\draw (.93,0) -- (1.65,0);
		
		\node[coupon,rotate=45]  at (-.72,.82)  {};
		\node[coupon,rotate=-45] at (-.72,-.82) {};
		\node[coupon,rotate=90]  at (.93,0)    {};
	\end{scope}
\end{tikzpicture}
	\caption{Replacing a vertex by joining up strands of the Jones-Wenzl idempotents}\label{fig:joneswenzl}
\end{figure}
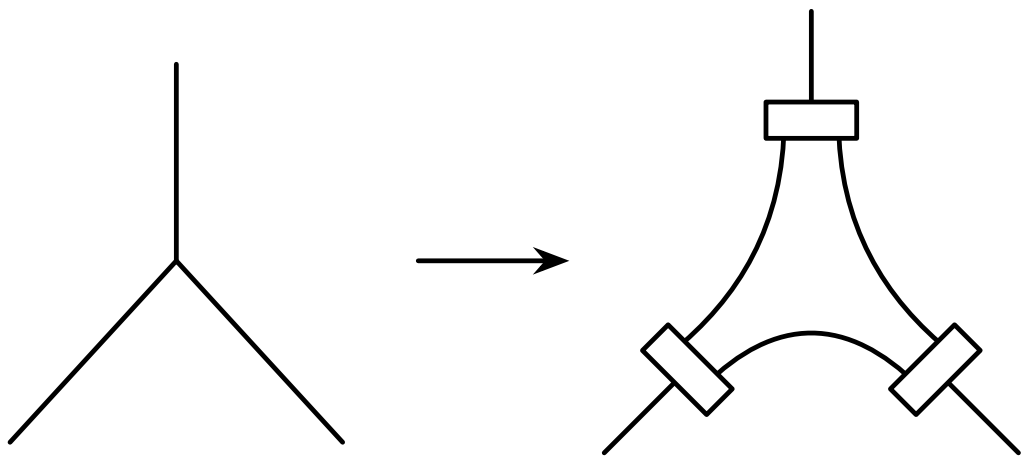

\begin{figure}
\begin{tikzpicture}[diagram]
	\draw (-2.4,1.2) -- (2.4,1.2)
	arc[start angle=90,end angle=-90,x radius=1.1,y radius=1.2]
	-- (-2.4,-1.2)
	arc[start angle=270,end angle=90,x radius=1.1,y radius=1.2]
	-- cycle;
	
	\foreach \x in {-1.05,1.05}{
		\begin{scope}[xshift=\x cm]
			\draw[meridian,hidden] (0,1.2)
			.. controls (-.22,.65) and (-.22,-.65) .. (0,-1.2);
		\end{scope}
	}
	
	\foreach \x in {-2.1,0,2.1}{
		\begin{scope}[xshift=\x cm]
			\draw[meridian,hidden] (0,-.11)
			.. controls (-.22,-.40) and (-.22,-1.10) .. (.22,-1.2);
		\end{scope}
	}
	
	\draw[meridian,hidden] (0,.064375)
	.. controls (-.22,.40) and (-.22,1.10) .. (.22,1.2);
	
	\draw[spine] (-1.28,0) -- (-.82,0);
	\draw[spine] (.82,0) -- (1.28,0);
	\foreach \x in {-2.1,0,2.1}{
		\draw[spine] (\x,0) ellipse [x radius=.82,y radius=.67];
	}
	
	\foreach \x in {-1.05,1.05}{
		\begin{scope}[xshift=\x cm]
			\draw[meridian] (0,1.2)
			.. controls (.22,.65) and (.22,-.65) .. (0,-1.2);
		\end{scope}
	}
	
	\draw[meridian] (0,.064375)
	.. controls (.28,.35) and (.32,.85) .. (.22,1.2);
	
	\foreach \x in {-2.1,0,2.1}{
		\begin{scope}[xshift=\x cm]
			\draw[meridian] (0,-.11)
			.. controls (.28,-.35) and (.32,-.85) .. (.22,-1.2);
			
			\draw (-.43,.16)
			.. controls (-.38,-.20) and (.38,-.20) .. (.43,.16);
			\draw (-.2816,-.0425)
			.. controls (-.14,.10) and (.14,.10) .. (.2816,-.0425);
		\end{scope}
	}
\end{tikzpicture}

\caption{A framed spine (in red, with the blackboard framing) for the handlebody $H$, and how it determines a disk system (in blue).}\label{fig:spine}
\end{figure}
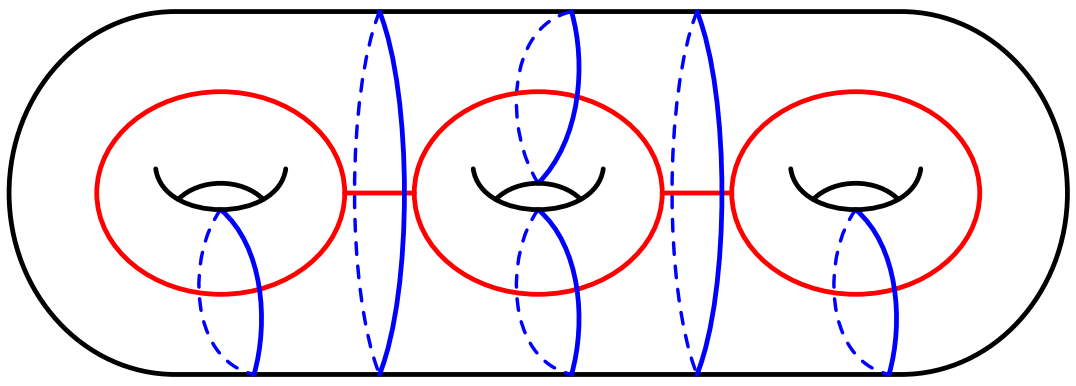
Throughout the paper, $R$ will be either $\Q(A)$ or $\C$ with a choice of $A\in \C^*$ not a root of unity. We also fix a closed, oriented $3$-manifold $M$ and a Heegaard splitting of $M$ of genus $g\geq 2$; we will write $M= H'\cup_{\partial H} H$. We also fix $\Gamma$ and $\Gamma'$ framed spines of $H$ and $H'$ respectively; for an admissible coloring $c$ of $\Gamma$, we will denote by $\phi_c$ the skein element in $S(H,R)$ given by replacing each edge $e$ of $\Gamma$ by the $c(e)$-th Jones-Wenzl idempotent and each vertex of $\Gamma$ as in Figure \ref{fig:joneswenzl}. We refer to the empty skein in $H$ by $\phi_0$, since coloring $\Gamma$ by the $0$ color on each edge gives the empty skein. Notice that because of our choice of $R$, all Jones-Wenzl idempotents are well defined (see for example \cite[Chapter 13]{Lic97}). We will use $\psi_c$ to refer to the skein element of $S(H',R)$ obtained by the same construction on $\Gamma'$.

Given a closed $3$-manifold $M$, we call the choice of a Heegaard splitting of $M$ together with a framed trivalent spine in each handlebody a \emph{Heegaard datum} for $M$, which we will denote by $\mathcal{H}$; for simplicity we will assume that the two spines are isomorphic as graphs.

The two spines $\Gamma$ and $\Gamma'$ provide disk systems for $H$ and $H'$, each disk corresponding to an edge of the graph (see Figure \ref{fig:spine}); we denote by $\alpha_e$ the curve in $\partial H$ bounding the disk corresponding to the edge $e$ of $\Gamma$, and by $\alpha'_e$ the analogous curve from $\Gamma'$. Furthermore, we choose $\Lambda\subseteq \partial H$ a parallel copy of $\Gamma$; for example, by identifying $H$ with a tubular neighborhood of $\Gamma$ and using the framing of $\Gamma$ to get a parallel in $\partial H$. Denote by $i(\gamma_1,\gamma_2)$ the geometric intersection count of the curves $\gamma_1,\gamma_2$. Similarly, we denote by $i(\gamma,\Lambda)$ the geometric intersection number of the graph $\Lambda$ with the curve $\gamma$, defined by isotoping $\gamma$ so that it is in general position with respect to $\Lambda$ and minimizing over all possible such choices.

\begin{definition}
	Given the Heegaard datum $\mathcal{H}$, we define $s(\mathcal{H}):=\textrm{max}_{e}i(\alpha'_e,\Lambda)$ and $w(\mathcal{H}):=\sum_{e,e'}i(\alpha_e,\alpha'_{e'})$.
	The \emph{intersection count} of $\mathcal{H}$ is defined by $B(\mathcal{H}):=2w(\mathcal{H})+s(\mathcal{H})$.
\end{definition}

From now on, we fix a Heegaard datum for $M$, and thus we omit $\mathcal{H}$ from the notation; we will simply use $s,$ $w$ and $B$ to refer to the quantities just introduced.

We need the following two, closely related, concepts of complexity on the skein modules at hand.

\begin{definition}\label{dfn:filtration}
	For a pants decomposition $\overline{\lambda}=\{\lambda_1,\dots,\lambda_{3g-3}\}$ of $\partial H$, we introduce the following filtration of $S(H,R)$:
	
	$$\mathcal{F}_N S(H,R)=\mathrm{Span}_{R}\left(\gamma\ | \ \gamma \subset \partial H \ \textrm{multicurve s.t. } \ \forall 1\leq i\leq 3g-3, i(\gamma,\lambda_i)\leq N \right)$$
	
	where curves in $\partial H$ are considered as curves in $H$ via the inclusion map.
	
	We define the same filtration for $S(\partial H,R)$. Notice that the filtration depends on the pants decomposition chosen; however we will only ever use the pants decomposition given by the spine of $H'$, and thus we will omit it from the notation.
	
\end{definition}

\begin{definition}\label{dfn:complexity}
	Let $H$ be a handlebody, $\Gamma$ a framed spine of $H$ and $\alpha_1,\dots,\alpha_{3g-3}$ the curves in $\partial H$ making up the corresponding disk system of $H$, with the convention that $\alpha_i$ is the dual curve to the edge $e_i\subseteq \Gamma$. Then for any admissible coloring $c$ of $\Gamma$ we denote by $\lvert c\rvert=\lvert \phi_c \rvert=\max_{i=1}^{3g-3}c(e_i)$ the \emph{degree} of $\phi_c$. 
\end{definition}
\begin{remark}
	Definitions \ref{dfn:filtration} and \ref{dfn:complexity} are very similar when applied to the same pants decomposition; if $\lvert \phi_c\rvert=N$ then $\phi_c\in \mathcal{F}_NS(H,R)$, and vice versa if $x\in \mathcal{F}_NS(H,R)$ then $x$ can be written as a linear combination of elements with $\lvert\phi_c\rvert\leq N$. The reason to give these two as separate definitions is that we have handlebodies $H$ and $H'$, each with a disk system; we will only consider the filtration $\mathcal{F}_N$ corresponding to the disk system of $H'$, and the degree with respect to the disk system of $H$, with the hope that this will be clearer.
\end{remark}
\section{Proof of the bound}

There are two main components in the proof of Theorem \ref{thm:upperbound}: one is a precise estimate of the dimension of $\mathcal{F}_N S(H',R)$ (Lemma \ref{lem:admissiblecounts}), and the other is a vacuum operator, whose existence is proven in Lemma \ref{lem:vacuum}. We refer to \cite[Chapter 9]{KL94} for an overview of the diagrammatic calculus used here.

We start by pointing out the following easy lemma.

\begin{lemma}\label{lem:shift}
	Suppose $x\in \mathcal{F}_M S(\partial H,R)$ and $\psi\in \mathcal{F}_N S(H',R)$; then $\psi\cdot x\in \mathcal{F}_{M+N}S(H',R)$.
\end{lemma}
\begin{proof}
	This follows immediately from the definitions of the filtration.
\end{proof}

\begin{lemma}\label{lem:admissiblecounts}
	For a handlebody $H'$ of genus $g$ with spine $\Gamma'$, there is a constant $K_{\Gamma'}>0$ such that $\dim \mathcal{F}_N S(H',R)=K_{\Gamma'}N^{3g-3}+O(N^{3g-4})$.
\end{lemma}
\begin{proof}
	Consider the skein elements $\psi_c$ given by admissibly coloring a spine $\Gamma'$ of $H'$; the ones with each color at most $N$ are linearly independent and span $\mathcal{F}_NS(H',R)$, so we just need to find an estimate of the number of admissible colorings with these properties. Concretely, we need to find the number of $(3g-3)$-tuples $(n_1,\dots,n_{3g-3})$ satisfying:
	\begin{itemize}
		\item $0\leq n_i\leq N$ for every $i$;
		\item $n_i,n_j,n_k$ satisfy all $3$ possible triangular inequalities whenever $i,j,k$ share a vertex;
		\item $n_i+n_j+n_k$ is even whenever $i,j,k$ share a vertex.
	\end{itemize}
	
	For now we ignore the last condition, and try to obtain a count of the tuples satisfying the first two; call this number $h_N$.
	
	Let $P$ be the polytope defined by:
		\begin{itemize}
		\item $0\leq x_i\leq 1$ for every $i$;
		\item $x_i+x_j\geq x_k$ whenever $i,j,k$ share a vertex.
	\end{itemize}
	
	Then $h_N$ is equal to the number of integral points of $NP$; by the theory of Ehrhart quasipolynomials (see for example \cite[Chapter 3]{BR07}) this is equal to $K'_{\Gamma'} N^{3g-3}+O(N^{3g-4})$. Adding an independent parity condition at a vertex cuts the leading term of the count in half; therefore, $\textrm{dim}\mathcal{F}_N S(H',R)=2^{-v}K'_{\Gamma'} N^{3g-3}+O(N^{3g-4})$ where $v$ is the number of independent parity conditions (this is equal to $2g-3$ as will be shown in the proof of Lemma \ref{lem:lowerbound}).
\end{proof}

\begin{lemma}\label{lem:lowerbound}
$K_{\Gamma'}\geq 2^{6-5g}$. 
\end{lemma}
\begin{proof}
	If a coloring assigns to each edge a color between $\lceil\frac{N}{2}\rceil$ and $N$, it will satisfy the triangular inequalities; there are $\left(\lfloor\frac{N}{2}\rfloor+1\right)^{3g-3}$ such colorings. 
	
	Each independent parity condition cuts this number in half, up to an error term in $O(N^{3g-4})$; by Euler characteristic arguments, a trivalent graph with $3g-3$ edges has $2g-2$ vertices, but one parity condition is superfluous, since each edge contributes twice to the sum of all vertex equations (so that this is always automatically even). Taking a spanning tree (which has one fewer edge than vertices) shows that the $2g-3$ remaining equations are independent, therefore there are $2^{6-5g}N^{3g-3}+O(N^{3g-4})$ admissible colorings with values between $\lceil \frac{N}{2}\rceil$ and $N$. 	
\end{proof}
\begin{remark}
	The contributions of $K_{\Gamma'}$ to the bound in Theorem \ref{thm:main} cancel out, but they will not for Theorem \ref{thm:boundary} addressing manifolds with boundary. We also note that Lemma \ref{lem:lowerbound} is a really crude estimate; a sharper one could be obtained by using the Verlinde formula, but this will both complicate the proof and make the upper bound of Theorem \ref{thm:boundary} more opaque (albeit, sharper).
\end{remark}

\begin{lemma}\label{lem:projection}
For every $N$ and every admissible coloring $c$ of the spine $\Gamma\subseteq H$ with $\lvert c\rvert \leq N$, there is an operator $P_{N,c}$ in $\mathcal{F}_{wN}S(\partial H,R)$ such that $P_{N,c}\phi_{c'}=\delta_{c,c'}\phi_c$ for all $c'$ with $\lvert c'\rvert\leq N$. 
\end{lemma}
\begin{proof}
	The key observation is that if $\alpha_e$ is the meridional curve in $\partial H$ dual to the edge $e$, then $\left(\alpha_e-\lambda_n\right)\phi_c=0$ if and only if $c(e)=n$, where $\lambda_n:=-A^{2n+2}-A^{-2n-2}$. This means that for a fixed $0\leq m\leq N$, the operator $T_{m,e}:=\prod_{0\leq n\leq N, n\neq m} \left(\alpha_e-\lambda_n\right)$ is such that if $\lvert c\rvert \leq N$, $T_{m,e}\phi_c=0$ if and only if $c(e)\neq m$; if $c(e)=m$ then $T_{m,e}\phi_c=\prod_{0\leq n\leq N, n\neq m}\left(\lambda_m-\lambda_n\right)\phi_c$ and the right-hand side coefficient is invertible in $R$ (since in particular, when $m\neq n$ we have $\lambda_n\neq \lambda_m$ in $R$). By construction, $T_{m,e}\in \mathcal{F}_{wN}$.
	
	We then define $\widetilde{P}_{N,c}$ to be $\prod_{e \textrm{ edge of }\Gamma}T_{c(e),e}$.
	
	 By the previous formula, $\widetilde{P}_{N,c}\phi_{c'}=0$ whenever $c\neq c'$ and $\lvert c'\rvert \leq N$, and $\widetilde{P}_{N,c}\phi_c=\eta_{c,N} \phi_c$ for some invertible element $\eta_{c,N}$. We then simply normalize $\widetilde{P}_{N,c}$ to obtain the desired operator $P_{N,c}$. 
	 
	 It remains to show that $P_{N,c}\in \mathcal{F}_{wN}S(\partial H,R)$. This follows directly from the definitions of $P_{N,c}$ and $w$; each $T_{m,e}$ is a polynomial in $\alpha_e$ of degree $N$, therefore each monomial in $P_{N,c}$ has total intersection number with the $\alpha'$ curves less than or equal to $wN$.

\end{proof}

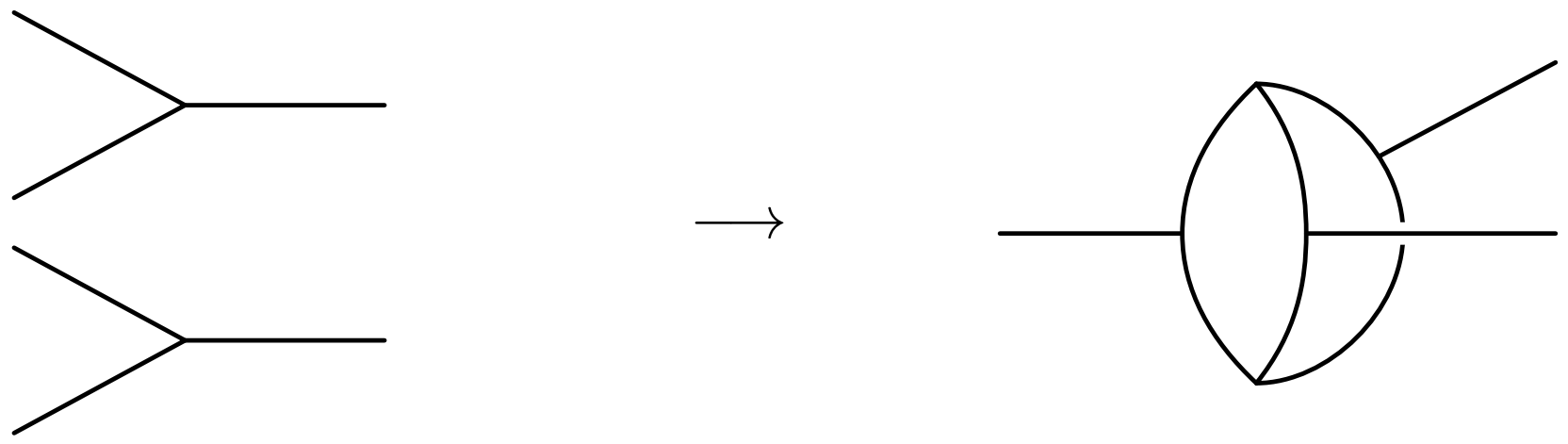
\begin{figure}
	\centering
	\begin{minipage}{.45\textwidth}
		\centering   
		\begin{tikzpicture}[diagram]
			\foreach \y in {0,-1.65}{
				\begin{scope}[yshift=\y cm]
					\draw (-1.2,.65) -- (0,0) -- (-1.2,-.65);
					\draw (0,0) -- (1.4,0);
				\end{scope}
			}
		\end{tikzpicture}
		
	\end{minipage} $\longrightarrow$  
	\begin{minipage}{.45\textwidth}
		\centering    	
		\begin{tikzpicture}[diagram]
			\coordinate (T) at (0,1.05);
			\coordinate (B) at (0,-1.05);
			\coordinate (L) at (-.52,0);
			\coordinate (C) at (.35,0);
			
			\draw (T) .. controls (-.32,.75) and (-.52,.4) .. (L)
			.. controls (-.52,-.4) and (-.32,-.75) .. (B);
			\draw (L) -- (-1.8,0);
			
			\draw (T) .. controls (.28,.7) and (.35,.35) .. (C)
			.. controls (.35,-.35) and (.28,-.7) .. (B);
			
			\draw (T) .. controls (.5,1.05) and (1.03,.55) .. (1.03,0)
			coordinate[pos=.65] (R)
			.. controls (1.03,-.55) and (.5,-1.05) .. (B);
			\draw (R) -- (2.1,1.2);
			
			\draw[white,line width=4.5pt] (.8,0) -- (1.25,0);
			\draw (C) -- (2.1,0);
		\end{tikzpicture}

	\end{minipage}
	\caption{Applying the fusion rule to every edge on the left side creates a bubble on the right.}\label{fig:bubble}
\end{figure}

To get the desired vacuum operator, we find an operator $M_c$ such that $M_c\cdot \phi_{c'}$ has a non-zero $\phi_0$-coefficient if and only if $c=c'$; we will then apply the operator of Lemma \ref{lem:projection} to obtain the desired operator.

\begin{lemma}\label{lem:vacuum}
	For every $N$ and every admissible coloring $c$ of the spine $\Gamma\subseteq H$ with $\lvert c\rvert \leq N$, there is an operator $V_{N,c}\in \mathcal{F}_{B N}S(\partial H,R)$ such that for every $c'$ with $\lvert c'\rvert \leq N$, $V_{N,c}\cdot \phi_{c'}=\delta_{c,c'}\phi_0$.
\end{lemma}
\begin{proof}
	Take $\phi_c$ and push it, via an isotopy of framed graphs, into the boundary of $H$ by identifying $H$ with a tubular neighborhood of $\Gamma$ and using the framing of $\Gamma$ to guide the isotopy; this is the same process used to define $\Lambda$ in the Heegaard datum. Call $M_c\in S(\partial H,R)$ the resulting skein element. By the definition of $s(\mathcal{H})$, $M_c\in \mathcal{F}_{sN}S(\partial H,R)$. 
	
	The first thing to notice is that in $S(H,R)$, $M_c\cdot \phi_{c'}$, written as a linear combination of the colored spines supported by $\Gamma$, contains as a summand $\phi_0$ with non-zero coefficient if and only if $c=c'$. First of all, apply the fusion rule to every corresponding pair of edges of $M_c$ and $\phi_{c'}$. When applying the fusion rule to the 3 corresponding pairs around a vertex $v$, the resulting graph will look like the bubble in Figure \ref{fig:bubble}; by elementary skein theory these bubbles can be replaced by a trivalent vertex at the cost of multiplying by some element of $R$; we do not make this coefficient explicit at first.
	
	After doing this for every pair of edges and after removing every bubble, we obtain $M_c\phi_{c'}=\sum_{\widetilde{c}} f_{\widetilde{c}} \phi_{\widetilde{c}}$, with the sum ranging over all admissible colorings $\widetilde{c}$ of $\Gamma$ such that for every edge $e$, $\lvert c(e)-c'(e)\rvert \leq \widetilde{c}(e)\leq c(e)+c'(e)$ (by the range of the fusion rule); notice that some coefficients could be equal to $0$.
	
	Therefore, if $c\neq c'$, we have written $M_c\phi_{c'}$ as a linear combination of colored spines not including $\phi_0$; since these form a basis, $M_c\phi_{c'}$ does not contain $\phi_0$. 
	
	Now suppose $c=c'$; we want to show that $M_c\cdot \phi_c$ contains $\phi_0$. In this case, the bubble coefficient in front of $\phi_0$ that corresponds to the vertex $v$ will simplify to be equal to $\theta(c(e_1),c(e_2),c(e_3))$ where $e_1,e_2,e_3$ are the edges around $v$ (possibly, with repetition). Since $c$ is admissible, this Theta coefficient is non-zero and a product and quotient of quantum integers, and thus invertible. Similarly, the coefficients arising from applying the fusion rule are also non-zero and invertible. Therefore, $\phi_0$ appears in $M_c\phi_c$ with non-zero coefficient. Furthermore, each coloring appearing in the fusion formula will have colors bounded above by $2N$; therefore, using the operator $P_{2N,0}$ from Lemma \ref{lem:projection}, we obtain that $P_{2N,0}M_c\cdot \phi_{c'}=\delta_{c,c'} \eta \phi_0$ for some invertible element $\eta$; normalizing it gives the operator $V_{N,c}$ required. To conclude, notice that $V_{N,c}\in \mathcal{F}_{BN}S(\partial H,R)$ by the property of the filtration.
\end{proof}

\begin{theorem}\label{thm:upperbound}
	Let $M$ be a closed manifold and $R$ be either $\Q(A)$ or $\C$ with $A$ not a root of unity. Suppose that $M=H\cup_{\partial H} H'$ is a Heegaard decomposition of genus $g\geq 2$, enhance this to a Heegaard datum $\mathcal{H}$ and denote by $B$ the quantity $s(\mathcal{H})+2w(\mathcal{H})$. Then, the dimension of $S(M,R)$ is less than or equal to $B^{3g-3}$.
\end{theorem}

\begin{proof}
	Suppose there are $d$ linearly independent elements of $S(M,R)$. Since $S(M,R)$ is spanned by the image of $S(\partial H,R)$, we pick $b_1,\dots,b_d\in S(\partial H,R)$ representing the $d$ linearly independent elements. We will let these elements act on the right on $S(H',R)$, and on the left on $S(H,R)$. Under these identifications, the elements $\psi_0\cdot b_i\cdot \phi_0$ are linearly independent elements of $S(M,R)$. Call $C_i$ the minimal integer such that $\psi_0\cdot b_i$ belongs to $\mathcal{F}_{C_i}S(H',R)$.
	
	 We claim that for any fixed $N$, the elements $y_{i,c,N}:=\psi_0\cdot b_i\cdot V_{N,c}$ are linearly independent in $S(H',R)$. To see this, suppose that some linear combination $\lambda=\sum_{i,c} f_{i,c}y_{i,c,N}=0$. Apply the pairing $S(H',R)\otimes S(H,R)\ra S(M,R)$ to $\lambda\otimes \phi_{c'}$, with $\lvert c'\rvert\leq N$. Expanding the sum gives $0=\sum_{i,c} f_{i,c}\psi_0\cdot b_i\cdot V_{N,c}\cdot \phi_{c'}$ and using Lemma \ref{lem:vacuum} gives $0=\sum_i f_{i,c'}\psi_0\cdot b_i\cdot \phi_0$ and linear independence implies that $f_{i,c'}=0$ for all $i$. Repeating this for every color $c'$ with $\lvert c'\rvert \leq N$ gives the desired result. 
	 
	 Let $h'(n)$ be the dimension of $\mathcal{F}_n S(H',R)$, and let $h(m)$ be the number of colorings of the spine $\Gamma$ with $\lvert c\rvert \leq m$. Notice that because $H$ and $H'$ are diffeomorphic and their spines are isomorphic, $h(n)=h'(n)$ for all $n$; we make the distinction in the hope of making it clear where elements belong. 
	 
	 By Lemma \ref{lem:admissiblecounts}, there is a constant $K_{\Gamma'}$ such that $h(N)=h'(N)=K_{\Gamma'}N^{3g-3}+O(N^{3g-4})$ for all $N$ big enough; therefore there are at least $d K_{\Gamma'} N^{3g-3}$ (up to an error term) linearly independent elements $y_{i,c,N}$ when $N$ is big enough.
	 
	 By Lemma \ref{lem:shift}, we have that $y_{i,c,N}$ belongs to $\mathcal{F}_{C_i+BN}S(H',R)$, so that every $y_{i,c,N}$ belongs to $\mathcal{F}_{C+BN}S(H',R)$ where $C=\max_iC_i$. However, $h(C+BN)=K_{\Gamma'} (BN+C)^{3g-3}+O(N^{3g-4})$, which means that $d K_{\Gamma'} N^{3g-3}\leq K_{\Gamma'} (BN+C)^{3g-3}+O(N^{3g-4})$ for all $N$ big enough, which implies that $d\leq B^{3g-3}$.
\end{proof}
\begin{remark}
	The arguments of the proof can be carried out, mutatis mutandis, for a genus $1$ Heegaard splitting, in which case the bound on the dimension would simply be equal to $B$; the main difference will be that for $g=1$, $h(N)=N+1$.
\end{remark}

\section{The case of manifolds with boundary}\label{sec:boundary}

For a manifold with boundary $M$, \cite{BD25} establishes that the skein module of $M$, localized at a pants decomposition, is finite dimensional. We recall here its definition.

\begin{definition}
	Let $M$ be a compact manifold with boundary $\partial M$, and let $\overline{\lambda}:=(\lambda_1,\dots,\lambda_n)$ be disjoint curves in $\partial M$. Then $R[\overline{\lambda}]$ forms a commutative subring of $S(\partial M,R)$; denote by $R(\overline{\lambda})$ the field of fractions of this subring. Then the \emph{localized skein module} of $M$ with respect to the curves $\lambda_1,\dots,\lambda_n$, denoted by $S_{\overline{\lambda}}(M)$, is the module $S(M,R)\otimes_{R[\overline{\lambda}]}R(\overline{\lambda})$. 
\end{definition}

All vector spaces (and their dimensions) from here on out are considered over $R(\overline{\lambda})$.

We adopt the slightly unusual convention that a pants decomposition of a surface is a maximal collection of disjoint, pairwise non-parallel, essential simple closed curves in the surface that cut it into spheres with holes (so a pants decomposition of a sphere is empty and of a torus is given by a single essential curve). Based on this, we can fill in every sphere boundary component on $M$ and assume that each component of $\partial M$ has genus at least $1$.

Now let $M$ be a manifold with boundary $\partial M$, and $\overline{\lambda}$ a pants decomposition of $\partial M$; it admits a Heegaard splitting $M=C\cup_{\partial H} H$; the compression body $C$ will have a negative boundary $\partial_-C$ identified with $\partial M$, and a positive boundary $\partial_+C$ identified with $\partial H$. We want to pick a graph in $C$ that will have the role of $\Gamma'$; to do so, we first use \cite[Lemma 3.2]{BD25} to extend $\overline{\lambda}$ to a pants decomposition of $\partial_+ C=\partial H$ in a compatible way (so that curves in the pants decomposition of $\partial_+ C$ either bound disks in $C$ or cobound annuli with curves in $\overline{\lambda}$) and then take $\Gamma'\subseteq \partial_+ C$ dual to the pants decomposition. We then fix $\Gamma$ a spine of $H$ with the same combinatorics as $\Gamma'$. Like before, we call this a Heegaard datum and we define the quantities $s,$ $w$ and $B$ exactly as before; the definition of the filtration on $S_{\overline{\lambda}}(C,R)$ is exactly the same as before, taking spans over $R(\overline{\lambda})$.

\begin{theorem}\label{thm:boundary}
	Take a compact $3$-manifold $M$ with a pants decomposition $\overline{\lambda}$ of $\partial M$, and a Heegaard datum $\mathcal{H}$ for $M$ with genus $g\geq 2$; then $\dim S_{\overline{\lambda}}(M,R)\leq (4B)^{3g-3}$.
\end{theorem}
To prove Theorem \ref{thm:boundary} we need a variation of Lemma \ref{lem:admissiblecounts} for compression bodies.
\begin{lemma}\label{lem:compressionbody}
	Let $C$ be a compression body, $\overline{\lambda}$ a pants decomposition of the negative boundary of $C$, and $S_{\overline{\lambda}}(C,R)$ be the localized skein module of $C$. Let $\overline{\lambda}'$ be a pants decomposition of the positive boundary of $C$ that is compatible with $\overline{\lambda}$. Then,
 $$\dim \mathcal{F}_N S_{\overline{\lambda}}(C,R)\leq 2^gN^{3g-3}+O(N^{3g-4})$$
 
 with $g\geq 2$ the genus of $\partial_+C$.
\end{lemma}
\begin{proof}
	It is shown in \cite[Lemma 3.3]{BD25} (technically, for $R=Q(A)$ but with a proof that works also for $R=C$, $A$ not a root of unity) that $\mathcal{F}_NS_{\overline{\lambda}}(C,R)$ is generated by multicurves whose Dehn-Thurston coordinates with respect to the pants decomposition $\overline{\lambda}'$ have intersection coordinates up to $N$ and twist coordinates  either $0$ or $1$. There are $2^{3g-3}(N+1)^{3g-3}$ possible Dehn-Thurston coordinates with each intersection at most $N$ and each twist between $0$ and $1$; however, the intersection coordinates must satisfy a parity condition at each pair of pants. There are $2g-3$ such independent conditions (using the same argument as Lemma \ref{lem:lowerbound}), therefore $\mathcal{F}_NS_{\overline{\lambda}}(C,R)$ has dimension at most $2^gN^{3g-3}+O(N^{3g-4})$.
\end{proof}

\begin{proof}[Proof of Theorem \ref{thm:boundary}]
	The argument of the proof follows exactly the same steps as the one of Theorem \ref{thm:main}, applying Lemma \ref{lem:compressionbody} to $C$ instead of Lemma \ref{lem:admissiblecounts}. The main difference is that in Theorem \ref{thm:upperbound}, the roles of the two handlebodies were the same, so we applied Lemma \ref{lem:admissiblecounts} to both and the contribution of $K_{\Gamma'}$ cancels out. Instead, here we apply Lemma \ref{lem:admissiblecounts} to estimate the dimension count of $H$ (with coefficients extended to $R(\overline{\lambda})$), and Lemma \ref{lem:compressionbody} for the dimension count of $C$. Therefore, we obtain that $\dim S_{\overline{\lambda}}(M,R)\leq \frac{2^g}{K_{\Gamma'}}B^{3g-3}$; applying Lemma \ref{lem:lowerbound} completes the proof. 
\end{proof}

\begin{remark}
	Corollary 1.7 of \cite{BD25} provides elements of the peripheral ideal of a manifold with boundary. This, in turn, gives recurrence relations for $J^K_n$, the colored Jones function of a framed knot $K$ (and, more in general, for the Reshetikhin-Turaev invariants and the colored Jones function of a framed knot in a closed manifold, see \cite{Bel25}). The degree of the recurrence relation (in the shifting variable) is bounded above by the dimension of $S_\lambda(E_K,\Q(A))$, where $E_K$ is the complement of $K$ and $\lambda$ is the longitude of $K$. Therefore, Theorem \ref{thm:boundary} provides an upper bound on the degree of a recurrence polynomial annihilating $J_n^K$.
\end{remark}

\printbibliography

\end{document}